




\input amstex
\documentstyle{amsppt}
\magnification\magstep1
\loadeusm

\def\a{\alpha}
\def\A{\eusm{A}}
\def\b{\beta}
\def\B{\eusm{B}(\eusm{H})}
\def\H{\eusm{H}}

\def\P{\Phi}
\def\R{\Bbb R}
\def\rng{\operatorname{rng}}
\def\S{\Psi}

\refstyle{A}
\widestnumber\key{ABCD}
\topmatter
\title
LINEAR MAPS ON FACTORS WHICH PRESERVE THE EXTREME POINTS OF THE UNIT
BALL
\endtitle
\rightheadtext{Linear maps which preserve the extreme points of the unit
ball}
\author
VANIA MASCIONI AND LAJOS MOLN\' AR
\endauthor
\thanks
This paper was written when the second author, holding a scholarship
of the Volkswagen-Stiftung, was a visitor at the
University of Paderborn, Germany. He is grateful to
Prof. K.-H.\linebreak Indlekofer for his kind hospitality.
The second author was partially supported also by the Hungarian National
Foundation for Scientific Research (OTKA), Grant No. T--016846
F--019322.
\endthanks
\address
Department of Mathematics,
The University of Texas at Austin,
Austin, TX 78712, U.S.A.
\endaddress
\email
{\tt mascioni\@math.utexas.edu}
\endemail
\address
Institute of Mathematics, Lajos Kossuth University, 4010
Debrecen, P.O.Box 12, Hungary
\endaddress
\email
{\tt molnarl\@math.klte.hu}
\endemail
\dedicatory
Dedicated to Zsuzsa \'Agnes Moln\'ar
\enddedicatory
\date
March 26, 1997
\enddate
\abstract
The aim of this paper is to characterize those linear maps
from a von Neumann factor $\A$ into itself which
preserve the extreme points of the unit ball of $\A$. For example, we
show that if $\A$ is infinite, then every such linear preserver can be
written as a fixed unitary operator times either
a unital *-homomorphism or a unital *-antihomomorphism.
\endabstract
\endtopmatter

\document
\head
Introduction and statements of the results
\endhead

Linear preserver problems deal with the question of
characterizing those linear maps on matrix algebras which leave
a certain subset, function or relation invariant. For example, in the
first mentioned case this means that the problem is to describe those
linear maps $\P$ on a matrix algebra $\eusm{M}$ for which
$\P(\eusm{S})\subset \eusm{S}$ holds true where $\eusm{S}$ is a given
subset of $\eusm{M}$. In fact, these problems
represent one of the most active research areas in matrix theory (see
the survey paper \cite{LiTs}). In the last decade considerable attention
has been paid to the infinite dimensional case as well, i.e. to linear
preserver problems concerning linear maps acting on operator algebras
rather than
matrix algebras (see the survey paper \cite{BrSe}). The linear preserver
problem we intend to investigate below is in an intimate connection
with the problem of unitary group preservers which are the linear maps
leaving the set of unitaries in $\eusm{M}$ invariant. The finite
dimensional case of this problem was
treated in \cite{Mar}, while the one on $\B$ (the
algebra of all bounded linear operators acting on the Hilbert space
$\H$) and on a general $C^*$-algebra were solved in \cite{Rai} and in
\cite{RuDy}, respectively. In \cite{Rai, Section 4} the problem of
characterizing those linear maps on $\B$ which preserve the extreme
points of the unit ball of $\B$ was implicitely raised and concerning
bijective
linear selfmaps of $\B$ which preserve the extreme points in question
in both directions (i.e. the maps as well as their inverses
are supposed to preserve the set of those extreme points) the author
obtained a complete description. The connection between
unitary group preservers on $\B$ and linear maps preserving the extreme
points of the unit ball of $\B$ is that in the first case our maps
preserve the set of all bijective partial isometries while in the
second case they preserve the set of all injective or surjective
partial isometries (see \cite{Hal, Sections 98, 99} and Lemma 1, Lemma
2 below).

The other motivation to our present investigations is the following.
In the paper \cite{LaMa} linear maps between $C^*$-algebras whose
adjoint preserve the extreme points of the dual ball were studied and
they turned out to give a valuable clue as to what objects may be
regarded as ``non-commutative composition operators''.
Now, it seems to be a natural problem to consider
linear maps in general, i.e. whithout the assumption of being
the adjoints of linear maps on a $C^*$-algebra, which preserve the
extreme points of the unit ball. For example, since $\B$ is the dual
space of the Banach algebra of all trace-class operators on $\H$ (which
is a highly non-$C^*$-algebra), in this particular but undoubtly
very important case the problem has nothing to do with
the one treated in \cite{LaMa}. In fact, as one can see below, we
follow a completely different approach to attack this problem.

As it is indicated in the abstract, we solve the problem of determing
all linear maps which preserve the extreme points of the unit ball in
the case when the underlying algebra is a von Neumann factor. One might
have the opinion that we should consider the problem for example
for general
$C^*$-algebras but an easy example shows that in that generality we
cannot expect almost anything. In fact, if the underlying algebra is the
$C^*$-algebra $\eusm{C}(\H)$ of all compact operators acting on an
infinite dimensional Hilbert space $\H$, then, since the unit ball of
this
Banach space has no extreme points at all, our preservers are
the linear maps on $\eusm{C}(\H)$ without any further properties.
So, to obtain a more satisfactory
result we have to suppose something more. In what follows we solve the
problem in the case when our preservers act on von Neumann factors
emphasizing the particular case of $\B$.

In view of the results as well as their proofs we have to remind the
definition of Jordan *-homomorphisms.
A linear map $J$ between *-algebras $\Cal A$ and $\Cal B$ is called
a Jordan *-homomorphism if
$$
\gather
J(x)^2=J(x^2)\\
J(x)^*=J(x^*)
\endgather
$$
hold true for every $x\in \Cal A$. Observe that by linearization,
i.e. replacing $x$ by $x+y$, the first equation above is equivalent to
$J(x)J(y)+J(y)J(x)=J(xy+yx)$ $(x,y\in \Cal A)$.

Let us now summarize the results of the paper.

\proclaim{Theorem 1}
Let $\A$ be an infinite factor. The linear map $\P:\A \to \A$ preserves
the extreme points of the unit ball of $\A$ if and only if either
there are a unitary operator $U\in \A$ and a unital *-homomorphism
$\S:\A \to \A$ such that $\P$ is of the form
$$
\P(A)=U\S(A) \qquad (A\in \A)
$$
or there are a unitary operator $U'\in \A$ and a unital
*-antihomomorphism $\S':\A \to \A$ such that $\P$ is of the form
$$
\P(A)=U'\S'(A) \qquad (A\in \A).
$$
\endproclaim

As a consequence of this result we immediately have the structure
of surjective linear selfmaps of $\B$ which preserve the extreme
points of the unit ball. In fact, Corollary 1 below is a significant
generalization
of a result of Rais \cite{Rai, Lemma 3 and Corollary 1} who obtained a
similar result but worked under the quite restrictive assumption
that the maps under consideration are bijective and
preserve the extreme points of the unit ball in both directions.

\proclaim{Corollary 1}
Let $\H$ be an infinite dimensional Hilbert space. Then the surjective
linear map $\P:\B \to \B$
preserves the extreme points of the unit ball of $\B$ if and only
if either there are unitaries $U,V\in \B$ such that $\P$ is of the form
$$
\P(A)=UAV \qquad (A\in \B)
$$
or there are antiunitaries $U',V'\in \B$ such that $\P$ is of the form
$$
\P(A)=U'A^*V' \qquad (A\in \B).
$$
\endproclaim

\remark{Remark}
We note that, as one can see from the proof of Corollary 1, it would
have been
sufficient to assume that the range of $\P$ contains a rank-one operator
instead of supposing that $\P$ is surjective.
\endremark

If the underlying Hilbert space is separable, then we can write
our linear preservers on $\B$ in a more detailed form than it was
obtained in the statement of Theorem 1.

\proclaim{Corollary 2}
Let $\H$ be a separable infinite dimensional Hilbert space. The linear
map $\P:\B\to \B$ preserves the extreme points of the
unit ball of $\B$ if and only if either there are a unitary operator $V$
and a collection $\{ U_\a\}$ of isometries with pairwise orthogonal
ranges which generate $\H$ such that $\P$ is
of the form
$$
\P(A)=V(\sum_\a U_\a AU^*_\a) \qquad (A\in \B)
$$
or there are a unitary operator $V$ and a family of anti-isometries $\{
V_\a\}$ with pairwise orthogonal ranges which generate $\H$ such that
$\P$ is of the form
$$
\P(A)=V(\sum_\b V_\b A^* V^*_\b) \qquad (A\in \B).
$$
\endproclaim

Our second main result describes our linear preservers in the case of
any finite von Neumann algebras.

\proclaim{Theorem 2}
Let $\A$ be a finite von Neumann algebra. The linear map $\P:\A \to \A$
preserves the extreme points of the unit ball of $\A$ if and only if
there exist a unitary operator $U\in \A$ and a unital Jordan
*-homomorphism $\S$ such that
$$
\P(A)=U\S(A) \qquad (A\in \A).
$$
\endproclaim

Concerning matrix algebras, we immediately have the last assertion of
the paper (cf. \cite{Mar}).

\proclaim{Corollary 3}
Let $\eusm{M}$ be the algebra of all complex $n\times n$ matrices. The
linear map $\P :\eusm{M}\to \eusm{M}$ preserves the extreme points of
the unit ball of $\eusm{M}$ if and only if there are unitary matrices
$U,V\in \eusm{M}$ such that $\P$ is either of the form
$$
\P(A)=UAV \qquad (A\in \eusm{M})
$$
or of the form
$$
\P(A)=UA^{tr}V \qquad (A\in \eusm{M})
$$
where ${}^{tr}$ denotes the transpose.
\endproclaim

\head
Proofs
\endhead

The statements of Lemma 1 and Lemma 2 are guessed to be well-known.
However, since
we have not found any trace of them in the bibliography of Kadison and
Ringrose \cite{KaRi1-2}, for the sake of completeness we present them
with proofs.

\proclaim{Lemma 1} Let $\A$ be a factor. The operator $A\in \A$
is an extreme point of the unit ball of $\A$ if and only if $A$
is either an isometry or a coisometry.
\endproclaim

\demo{Proof} It is well-known that in an arbitrary $C^*$-algebra
$\eusm{B}$, the extreme points of the unit ball are exactly those
partial isometries $W\in \eusm{B}$ for which
$(I-W^*W)\eusm{B}(I-WW^*)=\{ 0\}$ \cite{KaRi2, 7.3.1. Theorem}.
In a factor every two projections are comparable
\cite{KaRi2, 6.2.6. Proposition}.
Let, for example, $V\in \A$ be a partial isometry such that
$I-W^*W=V^*V$ and $VV^*$ is a subprojection of $I-WW^*$. Therefore,
we have $(V^*V)(V^*)(VV^*)=0$.
But $V$ is a partial isometry and hence $VV^*V=V$. Consequently, we
obtain that $0=(V^*V)(V^*VV^*)=V^*VV^*=V^*$ which implies $V=0$.
This gives us that $W$ is an isometry.

\qed
\enddemo

\demo{Proof of Theorem 1}
The sufficiency is trivial to check. Let us assume that
$\P$ preserves the extreme points of the unit ball of $\A$.
First observe that $\P$ is necessarily norm-continuous. Indeed,
since in an arbitrary $C^*$-algebra every self-adjoint operator
of norm $\leq 1$ is
the arithmetic mean of two unitaries, we easily obtain that $\|
\P\| \leq 2$.

Consider now the operator $V=\P(I)$. Since
$V$ is either an isometry or a coisometry, without loss of generality
we may and do suppose that $V^*V=I$.
Since the unitary group is arcwise connected in $\A$, we infer
that $\P(U)$ is an isometry for every unitary $U\in \A$.

Let us define a linear map $\S:\A
\to \A$ by $\S(A)=V^*\P(A)$ $(A\in \A)$.
In what follows we prove that $\S$ is a Jordan *-homomorphism.
To verify this, first observe that, by the fact that $\P$ sends
unitaries to isometries, we have
$$
\P(e^{itS})^*\P(e^{itS})=I \qquad (t\in \R)
$$
for every self-adjoint operator $S\in \A$.
Using the power series expansion of the exponential function as well
as its uniqueness, it is easy to conclude that
$$
\gather
\P(I)^*\P(S)-\P(S)^*\P(I)=0\\
-\frac{1}{2}\P(I)^*\P(S^2)+\P(S)^*\P(S)-\frac{1}{2}\P(S^2)^*\P(I)=0.
\endgather
$$
These identities imply that
$$
\P(S)^*\P(S)=\P(I)^*\P(S^2).
\tag 1
$$
Since $\S(S^2)=\P(I)^*\P(S^2)$, this shows that $\S$ preserves the
positive as well as the self-adjoint operators.
If we replace $S$ by $S+T$ in (1) where $S,T \in \A$ are self-adjoint,
then we obtain
$$
\P(S)^*\P(T)+\P(T)^*\P(S)=\P(I)^*\P(ST+TS)
$$
and one can check that this results in
$$
\P(A^*)^*\P(A)=\P(I)^*\P(A^2) \qquad (A\in \A).
$$
If we linearize this equation, i.e. replace $A$ by $A+B$,
we get
$$
\P(A^*)^*\P(B)+\P(B^*)^*\P(A)=\P(I)^*\P(AB+BA) \qquad (A,B\in \A).
\tag 2
$$
Let $P\in \A$ be an arbitrary projection. Since $\A$ is an infinite
factor, by \cite{StZs, E.4.11, p. 105} it follows that either $P\sim I$
or $I-P\sim I$. Suppose that this latter possibility is the case.
Then we have an isometry $U\in \A$ for which $UU^*=I-P$. By (2) we can
compute
$$
\gather
\P(U^*)^*\P(U^*)+\P(U)^*\P(U)=\P(I)^*\P(UU^*+U^*U)=\\
=\P(I)^*\P(2I-P)=2-\S(P). \tag 3
\endgather
$$
Since, by our assumption on $\P$, $\P(U^*)^*\P(U^*)$ and $\P(U)^*\P(U)$
are projections, we infer that there are projections $Q_1,Q_2\in \A$
such that $\S(P)=Q_1+Q_2$. On the other hand, by the positivity
preserving property of $\S$ and (1) we have
$$
\S(P)^2=\S(P)^*\S(P)=\P(P)^*VV^*\P(P)\leq \P(P)^*\P(P)=\S(P)
$$
which results in $Q_1Q_2+Q_2Q_1\leq 0$. Since $Q_1,Q_2$ are projections,
we easily obtain $Q_1Q_2=Q_2Q_1=0$ and this gives us that $\S(P)$ is a
projection.
It is now a standard argument to show that $\S$ is a Jordan
*-homomorphism.
Indeed, if $P_1,P_2\in \A$ are mutually orthogonal projections, then we
know that
$\S(P_1+P_2)=\S(P_1)+\S(P_2)$ is a projection, too. Therefore, $\S(P_1)$
and $\S(P_2)$ are also orthogonal, i.e. $\S$ preserves the
orthogonality between projections. If $P_1,\dots, P_n\in \A$ are
mutually
orthogonal projections and $\l_1,\dots,\l_n\in \R$, then we have
$$
\bigl[\S(\sum_{k=1}^n\l_kP_k)\bigr]^2=
\bigl[\sum_{k=1}^n\l_k\S(P_k)\bigr]^2=
\sum_{k=1}^n\l_k^2\S(P_k)=
\S((\sum_{k=1}^n\l_kP_k)^2).
$$
Since $\S$ is norm-continuous, by the spectral theorem we have
$$
\S(S)^2=\S(S^2)
$$
for every self-adjoint operator $S\in \A$. Using the trick of
linearization we readily obtain
$$
\S(A)^2=\S(A^2)  \qquad (A\in \A).
$$
Since $\S$ is linear and positivity preserving, we obviously have
$$
\S(A)^*=\S(A^*) \qquad (A\in \A).
$$
Consequently, we infer that $\S$ is a Jordan *-homomorphism.

Our next claim is that $\P(A)=V\S(A)$ $(A\in \A)$.
To this end, let $U\in \A$ be an arbitrary unitary operator.
Then $W=\P(U)$ is an isometry and we have
$$
\S(U)^*\S(U)+\S(U)\S(U)^*=\S(U^*U+UU^*)=\S(2I)=2I,
$$
i.e.
$$
W^*VV^*W+V^*WW^*V=2I.
\tag 4
$$
Since $W^*VV^*W\leq I$ and
$V^*WW^*V\leq I$, by (4) it follows that $W^*VV^*W= I$ and
$V^*WW^*V=I$. Since $VV^*,WW^*$ are projections, these relations
imply $\rng W\subset \rng V$ and $\rng V\subset \rng W$. Consequently,
$\rng \P(U)=\rng V$ for every unitary operator $U\in \A$. Since the
linear
span of the unitaries in $\A$ is $\A$, it follows that $\rng
\P(A)\subset \rng V$ for every $A\in \A$, and
this gives us that $V\S(A)=VV^*\P(A)=\P(A)$ $(A\in \A)$.

We now prove that $V$ is unitary.
Since $\S$ is a unital Jordan *-homomorphism, by \cite{Sto, Theorem 3.3}
there is a projection $E$ in the centre of the von Neumann algebra
generated by the image of $\S$, such that the maps
$\S_1$ and $\S_2$ on $\A$ defined by $\S_1(A) = \S(A)E$  and
$\S_2(A) = \S(A) (I-E)$  are a *-homomorphism and a
*-antihomomorphism, respectively.
Now, if we suppose on the contrary that
$\rng V\neq \H$, then by the extreme point preserving property of $\P$
and the equality $\P(A)=V\S(A)$ it follows that $\P$ sends every
isometry and coisometry to a proper isometry. Let
$W^*W=I$ or $WW^*=I$. In both cases we have
$I=\S(W)^*V^*V\S(W)=\S(W)^*\S(W)=\S_1(W^*W)+\S_2(WW^*)$.
This readily imples that
$$
\S_1(W^*W)=E\quad \text{ and }\quad \S_2(WW^*)=I-E.
\tag 5
$$
Let $P\in \A$ be a projection such that $P\sim I$ and $I-P\sim I$
\cite{KaRi2, 6.3.3. Lemma (Halving)}. Using (5) we have
$\S_1(P)=\S_1(I-P)=E,\, \S_2(P)=\S_2(I-P)=I-E$. Since the
operators appearing here are all projections, we deduce
$\S_1(P)=\S_1(I)=0$ and $\S_2(P)=\S_2(I)=0$.
It now follows that $\S_1=0$ and $\S_2=0$ which
is an obvious contradiction. Therefore, we obtain that $V$ is unitary.

It only remains to prove that either $\S_1=0$ or $\S_2=0$.
Just as before, let $P\in
\A$ be a projection for which $P\sim I$ and $I-P\sim I$. Let $W\in \A$
be such that $W^*W=I$ and $WW^*=P$. Then $\S(W)$ is either an isometry
or a coisometry. Suppose that this latter one is the case. Then we
have
$$
I=\S(W)\S(W)^*=\S_1(WW^*)+\S_2(W^*W)
$$
which gives us that $\S_1(P)=\S_1(I)$. If $W'$ is an isometry for which
$I-P=W'{W'}^*$, then we have $\S_1(W')\S_1(W')^*=0$. Therefore,
$\S_1(W')=0$ and hence $0=\S_1(W')^*\S_1(W')=\S_1(I)$. It follows that
$\S_1=0$ and thus we obtain that $\S$ is a *-antihomomorphism.
In case $\S(W)$ is an isometry, one can argue in a very similar way.

\qed
\enddemo

\demo{Proof of Corollary 1}
Using Theorem 1, without a serious loss of generality we may and do
suppose that our map $\P$ is a surjective unital *-homomorphism of $\B$.
We assert that $\P$ is injective.
Let $P$ be a rank-one projection. Then there is an operator $A\in \B$
such that $\P(A)=P$. Since $\P(A^*A)=\P(A)^*\P(A)=P$, we can assume that
our operator $A$ is positive. Let us consider its spectral resolution.
By the monotonicity and continuity of $\P$ it follows that there is a
Borel subset $B$ of the spectrum of $A$ having a positive distance from
0 for which $E(B)$ ($E$ is the spectral measure corresponding to $A$)
has nonzero image under $\P$. Since there is a positive constant $c$
for which $cE(B)\leq A$, we obtain $0\neq c\P(E(B))\leq P$. By the
minimality of $P$ we have $\P(E(B))=P$. Now, if $\P$ is not
injective, then its kernel, being a nontrivial ideal of $\B$, contains
the ideal of all finite rank operators. Hence, the projection $E(B)$ is
infinite dimensional. But in this case $E(B)$ is the sum of two
orthogonal projections $Q_1, Q_2$ which are both equivalent to $E(B)$.
Using the minimality of $P$ once again, we see that either $\P(Q_1)=0$
or
$\P(Q_2)=0$. Now, $Q_1$ and $Q_2$ are equivalent to each other and $\P$
is a *-homomorphism, thus if any of $Q_1,Q_2$ is in the kernel of $\P$,
then so is the other one. Hence, we deduce $0=\P(Q_1)+\P(Q_2)=P$ which
is a contradiction. Consequently, we obtain the injectivity of $\P$ and so
$\P$ is a *-automorphism of $\B$. To conclude, it is a folk result
that in this case there is a unitary operator $U\in \B$ such that
$\P(A)=UAU^*$ $(A\in \B)$. This completes the proof.

\qed
\enddemo

\demo{Proof of Corollary 2}
Using the notation of Theorem 1 we can suppose without loss of
generality that $\P:\B \to \B$ is a unital
*-homomorphism. Now, using the separability of $\H$ and the
classical result \cite{KaRi2,
10.4.14. Corollary} on the form of the representations of $\B$ on
$\B$, we obtain the assertion.

\qed
\enddemo

\proclaim{Lemma 2}
Let $\A$ be a finite von Neumann algebra. The extreme points of the unit
ball of $\A$ are exactly the unitaries in $\A$.
\endproclaim

\demo{Proof}
By \cite{KaRi2, 7.3.1. Theorem} every unitary operator in $\A$ is an
extreme point of the unit ball of $\A$. Referring to that theorem again,
suppose that
$V\in \A$ is a partial isometry for which $(I-V^*V)\A(I-VV^*)=\{ 0\}$.
Let $V^*V=E, VV^*=F$.
Since $\A$ is finite, by \cite{KaRi2, 6.9.6.}
we have $I-E\sim I-F$. Let $W\in \A$ be a partial isometry for which
$I-E=W^*W$ and $I-F=WW^*$. Since $W^*W\A WW^*=\{ 0\}$, it follows that
$0=W^*WW^*WW^*=W^*WW^*=W^*$. Consequently, we have $E=F=I$ which means
that $V$ is unitary.

\qed
\enddemo

\demo{Proof of Theorem 2}
If $\S$ is a unital Jordan *-homomorphism on $\A$, then taking the fact
that $\S$ is necessarily a contraction into consideration, by the
equality
$$
\S(U)^*\S(U)+\S(U)\S(U)^*=\S(2I)=2I
$$
we easily obtain that $\S(U)$ is unitary for every unitary operator
$U\in \A$. On the other hand, if $\P$ preserves the extreme
points of $\A$, i.e. the unitary group of the $C^*$-algebra $\A$, then
by \cite{RuDy, Corollary 2} we obtain the other part of the assertion.

\qed
\enddemo

\demo{Proof of Corollary 3}
The sufficiency is obvious. To the necessity we may suppose
that $\P$ is a unital Jordan *-homomorphism. Since
the algebra under consideration is finite dimensional, we obtain
that $\P$ is an injective and hence surjective Jordan *-homomorphism,
i.e. a Jordan *-automorphism. By a well-known theorem of Herstein
\cite{Her} $\P$ is either a *-automorphism or a *-antiautomorphism and,
just as in the proof of Corollary 1, we are done.

\qed
\enddemo

\enddocument

\end